\documentclass[10pt,twoside]{siamltex}
\usepackage{amsmath}
\usepackage{amsfonts,epsfig}

\newtheorem{algorithm}[theorem]{Algorithm}

\input{tcilatex}
\begin{document}

\title{ Condensation of Determinants}
\author{ Abdelmalek Salem\thanks{%
Department of mathematics, University Centre of Tebessa 12002 Algeria,
(a.salem@gawab.com)} \and Kouachi Said\thanks{%
University Centre of Khenchela 40100 Algeria, (kouachi.said@caramail.com)}}
\maketitle

\begin{abstract}
In this paper we tried to condense the determinant of $n$ square matrix to
the determinant of $(n-1)$ square matrix with the mathematical proof.
\end{abstract}





\pagestyle{myheadings} 
\markboth{S.\ Abdelmalek and S.\
Kouachi}{Condensation of Determinants}



\section{Introduction}

We can write the well-known algorithm of Dodgson, concerning the\ $n\ $%
square matrix $A=\left( a_{i,j}\right) _{1\leq i,j\leq n}$ as follows:

\begin{equation}
\det \left[ \left( a_{i,j}\right) _{1\leq i,j\leq n}\right] \det \left[
\left( a_{i,j}\right) _{\substack{ i\neq k,l  \\ j\neq k,l}}\right] =\det %
\left[ 
\begin{array}{cc}
\det \left[ \left( a_{i,j}\right) _{\substack{ i\neq l  \\ j\neq l}}\right]
\smallskip & \det \left[ \left( a_{i,j}\right) _{\substack{ i\neq l  \\ %
j\neq k}}\right] \\ 
\det \left[ \left( a_{i,j}\right) _{\substack{ i\neq k  \\ j\neq l}}\right]
& \det \left[ \left( a_{i,j}\right) _{\substack{ i\neq k  \\ j\neq k}}\right]%
\end{array}%
\right] ,  \label{(1)}
\end{equation}%
for all $k,l=1,...,n$ considering $k<l.$(see S.Kouachi, S.Abdelmalek and
B.Rebai \cite{Kouachi})

This formula enables us condense the determinant of $n$ square matrix to the
determinant of $2$ square matrix. The elements of $2$ square matrix are the
determinants of $(n-1)$ square matrix.

In the same way we try to create a formula that enables us condense the
determinant of $n$ square matrix to the determinant of $(n-1)$ square
matrix. The elements of $(n-1)$ square matrix are the determinants of $2$
square matrix.

For exampel $n=7$

\smallskip $\left( 2\right) ^{7-2}\left\vert 
\begin{array}{ccccccc}
2 & 5 & 4 & 7 & 6 & 1 & 2 \\ 
0 & 1 & 3 & 8 & 8 & 1 & 5 \\ 
9 & 4 & 7 & 8 & 9 & 8 & 6 \\ 
7 & 8 & 4 & \sqrt{3} & 2 & 0 & 8 \\ 
11 & 2 & 5 & 4 & 5 & \frac{1}{2} & 5 \\ 
5 & 7 & 8 & 6 & 1 & 0 & 5 \\ 
9 & 2 & 3 & 5 & 8 & 5 & 3%
\end{array}%
\right\vert =$

$\left\vert 
\begin{array}{cccccc}
\left\vert 
\begin{array}{cc}
2 & 5 \\ 
0 & 1%
\end{array}%
\right\vert \smallskip & \left\vert 
\begin{array}{cc}
2 & 4 \\ 
0 & 3%
\end{array}%
\right\vert & \left\vert 
\begin{array}{cc}
2 & 7 \\ 
0 & 8%
\end{array}%
\right\vert & \left\vert 
\begin{array}{cc}
2 & 6 \\ 
0 & 8%
\end{array}%
\right\vert & \left\vert 
\begin{array}{cc}
2 & 1 \\ 
0 & 1%
\end{array}%
\right\vert & \left\vert 
\begin{array}{cc}
2 & 2 \\ 
0 & 5%
\end{array}%
\right\vert \\ 
\left\vert 
\begin{array}{cc}
2 & 5 \\ 
9 & 4%
\end{array}%
\right\vert \smallskip & \left\vert 
\begin{array}{cc}
2 & 4 \\ 
9 & 7%
\end{array}%
\right\vert & \left\vert 
\begin{array}{cc}
2 & 7 \\ 
9 & 8%
\end{array}%
\right\vert & \left\vert 
\begin{array}{cc}
2 & 6 \\ 
9 & 9%
\end{array}%
\right\vert & \left\vert 
\begin{array}{cc}
2 & 1 \\ 
9 & 8%
\end{array}%
\right\vert & \left\vert 
\begin{array}{cc}
2 & 2 \\ 
9 & 6%
\end{array}%
\right\vert \\ 
\left\vert 
\begin{array}{cc}
2 & 5 \\ 
7 & 8%
\end{array}%
\right\vert \smallskip & \left\vert 
\begin{array}{cc}
2 & 4 \\ 
7 & 4%
\end{array}%
\right\vert & \left\vert 
\begin{array}{cc}
2 & 7 \\ 
7 & \sqrt{3}%
\end{array}%
\right\vert & \left\vert 
\begin{array}{cc}
2 & 6 \\ 
7 & 2%
\end{array}%
\right\vert & \left\vert 
\begin{array}{cc}
2 & 1 \\ 
7 & 0%
\end{array}%
\right\vert & \left\vert 
\begin{array}{cc}
2 & 2 \\ 
7 & 8%
\end{array}%
\right\vert \\ 
\left\vert 
\begin{array}{cc}
2 & 5 \\ 
11 & 2%
\end{array}%
\right\vert \smallskip & \left\vert 
\begin{array}{cc}
2 & 4 \\ 
11 & 5%
\end{array}%
\right\vert & \left\vert 
\begin{array}{cc}
2 & 7 \\ 
11 & 4%
\end{array}%
\right\vert & \left\vert 
\begin{array}{cc}
2 & 6 \\ 
11 & 5%
\end{array}%
\right\vert & \left\vert 
\begin{array}{cc}
2 & 1 \\ 
11 & \frac{1}{2}%
\end{array}%
\right\vert & \left\vert 
\begin{array}{cc}
2 & 2 \\ 
11 & 5%
\end{array}%
\right\vert \\ 
\left\vert 
\begin{array}{cc}
2 & 5 \\ 
5 & 7%
\end{array}%
\right\vert \smallskip & \left\vert 
\begin{array}{cc}
2 & 4 \\ 
5 & 8%
\end{array}%
\right\vert & \left\vert 
\begin{array}{cc}
2 & 7 \\ 
5 & 6%
\end{array}%
\right\vert & \left\vert 
\begin{array}{cc}
2 & 6 \\ 
5 & 1%
\end{array}%
\right\vert & \left\vert 
\begin{array}{cc}
2 & 1 \\ 
5 & 0%
\end{array}%
\right\vert & \left\vert 
\begin{array}{cc}
2 & 2 \\ 
5 & 5%
\end{array}%
\right\vert \\ 
\left\vert 
\begin{array}{cc}
2 & 5 \\ 
9 & 2%
\end{array}%
\right\vert & \left\vert 
\begin{array}{cc}
2 & 4 \\ 
9 & 3%
\end{array}%
\right\vert & \left\vert 
\begin{array}{cc}
2 & 7 \\ 
9 & 5%
\end{array}%
\right\vert & \left\vert 
\begin{array}{cc}
2 & 6 \\ 
9 & 8%
\end{array}%
\right\vert & \left\vert 
\begin{array}{cc}
2 & 1 \\ 
9 & 5%
\end{array}%
\right\vert & \left\vert 
\begin{array}{cc}
2 & 2 \\ 
9 & 3%
\end{array}%
\right\vert%
\end{array}%
\right\vert .$

Finally:

$32\left\vert 
\begin{array}{ccccccc}
2 & 5 & 4 & 7 & 6 & 1 & 2 \\ 
0 & 1 & 3 & 8 & 8 & 1 & 5 \\ 
9 & 4 & 7 & 8 & 9 & 8 & 6 \\ 
7 & 8 & 4 & \sqrt{3} & 2 & 0 & 8 \\ 
11 & 2 & 5 & 4 & 5 & \frac{1}{2} & 5 \\ 
5 & 7 & 8 & 6 & 1 & 0 & 5 \\ 
9 & 2 & 3 & 5 & 8 & 5 & 3%
\end{array}%
\right\vert =\left\vert 
\begin{array}{cccccc}
2 & 6 & 16 & 16 & 2 & 10 \\ 
-37 & -22 & -47 & -36 & 7 & -6 \\ 
-19 & -20 & 2\sqrt{3}-49 & -38 & -7 & 2 \\ 
-51 & -34 & -69 & -56 & -10 & -12 \\ 
-11 & -4 & -23 & -28 & -5 & 0 \\ 
-41 & -30 & -53 & -38 & 1 & -12%
\end{array}%
\right\vert .$

For this purpose, we need some notations:

The $\left( n-k\right) \times \left( n-l\right) $ matrix obtained from $A$
by removing the $i_{1}^{th},i_{2}^{th}...i_{k}^{th}$ rows and the $%
j_{1}^{th},$ $j_{2}^{th};...j_{l}^{th}$ columns is denoted by $\left(
a_{i,j}\right) _{\substack{ i\neq i_{1},i_{2},...i_{k}  \\ j\neq
j_{1},j_{2},...j_{l}}}$.\newline
We denote by $\underset{\alpha \leq i,j\leq \beta }{\det }\left[ \left(
a_{i,j}\right) _{\substack{ i\neq i_{1},i_{2},...i_{k}  \\ j\neq
j_{1},j_{2},...j_{k}}}\right] $, $k\leq \beta -\alpha $ to the determinant
of the $\left( \beta -\alpha -k+1\right) $ square matrix\ obtained from $A$
by removing the$\left( \alpha -1\right) $ first rows and columns, by
removing the $\left( n-\beta \right) $ last rows and columns and by removing 
$i_{1}^{th},i_{2}^{th}...i_{k}^{th}$ rows and the $j_{1}^{th},$ $%
j_{2}^{th};...j_{k}^{th}$ columns.

\section{RESULTS}

We need a lemma.

\begin{lemma}
\label{lemma1}If $a_{11}=0$, thus we get the following formula: 
\begin{equation}
\underset{1\leq i,j\leq n-1}{\det }\left[ \det \left[ 
\begin{array}{cc}
a_{1,1} & a_{1,j+1} \\ 
a_{i+1,1} & a_{i+1,j+1}%
\end{array}%
\right] \right] =0\text{.}  \label{(2.1)}
\end{equation}
\end{lemma}

\begin{proof}
To prove the formula (\textrm{\ref{(2.1)})}, we let $a_{11}=0$. Thus, the
first term of this formula will be as follows:

$\underset{1\leq i,j\leq n-1}{\det }\left[ \det \left[ 
\begin{array}{cc}
0 & a_{1,j+1} \\ 
a_{i+1,1} & a_{i+1,j+1}%
\end{array}%
\right] \right] =\underset{1\leq i,j\leq n-1}{\det }\left[ \left(
-a_{1,j+1}a_{i+1,1}\right) \right] =$

$\left( -1\right) ^{n-1}\left( \underset{i=2}{\overset{n}{\tprod }}%
a_{1,i}\right) \left( \underset{i=2}{\overset{n}{\tprod }}a_{i,1}\right) 
\overset{n-1}{\overbrace{\left\vert 
\begin{array}{ccccc}
1 & 1 & ... & 1 & 1 \\ 
1 & 1 & ... & 1 & 1 \\ 
\vdots & \vdots & \ddots & \vdots & \vdots \\ 
1 & 1 & ... & 1 & 1 \\ 
1 & 1 & ... & 1 & 1%
\end{array}%
\right\vert }}.$

And as it known that $\left\vert 
\begin{array}{ccccc}
1 & 1 & ... & 1 & 1 \\ 
1 & 1 & ... & 1 & 1 \\ 
\vdots & \vdots & \ddots & \vdots & \vdots \\ 
1 & 1 & ... & 1 & 1 \\ 
1 & 1 & ... & 1 & 1%
\end{array}%
\right\vert =0$, and like this we have finished the proof of lemma (\textrm{%
\ref{lemma1}).}
\end{proof}

One of the main results of the paper is the following:

\begin{theorem}
\label{thor1}Let the $n$ square matrix $A=\left( a_{\left( i,j\right)
}\right) _{1\leq i,j\leq n}.$\newline
For all $n>2$ , the following formula is realised\newline
\begin{equation}
\left( a_{1,1}\right) ^{n-2}\underset{1\leq i,j\leq n}{\det }\left[ \left(
a_{i,j}\right) _{1\leq i,j\leq n}\right] =\underset{1\leq i,j\leq n-1}{\det }%
\left[ \det \left[ 
\begin{array}{cc}
a_{1,1} & a_{1,j+1} \\ 
a_{i+1,1} & a_{i+1,j+1}%
\end{array}%
\right] \right] .  \label{(2.2)}
\end{equation}%
\newline
We notice that this formula enables us condense the determinant of $n$
square matrix to the determinant of $\left( n-1\right) $ square matrix. The
elements of $\left( n-1\right) $ square matrix are the determinants of $2$
square matrix.
\end{theorem}

\begin{proof}
To prove formula (\textrm{\ref{(2.2)})} there are two cases:

The first case when $a_{1,1}=0$, the proof of formula (\textrm{\ref{(2.2)})}
is the same as the one of lemma (\textrm{\ref{lemma1})}.

The second case when $a_{1,1}\neq 0$, we prove the formula (\textrm{\ref%
{(2.2)})} inductively.

For $n=3$, we find that the proof of the formula (\textrm{\ref{(2.2)})} is
evident.

For $n=4$, we find that the proof of the formula (\textrm{\ref{(2.2)})} is
evident with simple calculations.

When $n>4$, we suppose the formula (\textrm{\ref{(2.2)})} is correct for $%
\left( n-1\right) $ and we prove it for $n.$

To prove the case of $\left( n>4\right) $ we use formula (\textrm{\ref{(1)})}
, and to choose $\det \left[ \left( a_{i,j}\right) _{\substack{ i\neq k,l 
\\ j\neq k,l}}\right] \neq 0$ with $k>1$ (this choice is possible).

We assume without loss of generality that $\det \left[ \left( a_{i,j}\right) 
_{\substack{ i\neq n-1,n  \\ j\neq n-1,n}}\right] \neq 0$, and this means
that $k=n-1,$ $l=n$, thus formula \textrm{\ref{(1)}}, will be as follows:

\begin{equation}
\det \left[ \left( a_{i,j}\right) _{1\leq i,j\leq n}\right] \det \left[
\left( a_{i,j}\right) _{\substack{ i\neq n-1,n  \\ j\neq n-1,n}}\right]
=\det \left[ 
\begin{array}{cc}
\det \left[ \left( a_{i,j}\right) _{\substack{ i\neq n  \\ j\neq n}}\right]
& \det \left[ \left( a_{i,j}\right) _{\substack{ i\neq n  \\ j\neq n-1}}%
\right] \\ 
\det \left[ \left( a_{i,j}\right) _{\substack{ i\neq n-1  \\ j\neq n}}\right]
& \det \left[ \left( a_{i,j}\right) _{\substack{ i\neq n-1  \\ j\neq n-1}}%
\right]%
\end{array}%
\right]  \label{(2.3)}
\end{equation}

We apply the formula (\textrm{\ref{(2.2)})} for $\left( n-1\right) $ on the
second side of formula (\textrm{\ref{(2.3)})}, so we get:

\begin{equation}
\left( a_{11}\right) ^{n-3}\det \left[ \left( a_{i,j}\right) _{\substack{ %
i\neq n  \\ j\neq n}}\right] =\underset{1\leq i,j\leq n-2}{\det }\left[ \det %
\left[ 
\begin{array}{cc}
a_{1,1} & a_{1,j+1} \\ 
a_{i+1,1} & a_{i+1,j+1}%
\end{array}%
\right] \right]  \label{(2.4)}
\end{equation}

\begin{equation}
\left( a_{11}\right) ^{n-3}\det \left[ \left( a_{i,j}\right) _{\substack{ %
i\neq n-1  \\ j\neq n-1}}\right] =\underset{1\leq i,j\leq n-1}{\det }\left[
\det \left[ 
\begin{array}{cc}
a_{1,1} & a_{1,j+1} \\ 
a_{i+1,1} & a_{i+1,j+1}%
\end{array}%
\right] \right] _{\substack{ i\neq n-2  \\ j\neq n-2}}  \label{(2.5)}
\end{equation}

\begin{equation}
\left( a_{11}\right) ^{n-3}\det \left[ \left( a_{i,j}\right) _{\substack{ %
i\neq n  \\ j\neq n-1}}\right] =\underset{1\leq i,j\leq n-1}{\det }\left[
\det \left[ 
\begin{array}{cc}
a_{1,1} & a_{1,j+1} \\ 
a_{i+1,1} & a_{i+1,j+1}%
\end{array}%
\right] \right] _{\substack{ i\neq n-1  \\ j\neq n-2}}  \label{(2.6)}
\end{equation}

\begin{equation}
\left( a_{11}\right) ^{n-3}\det \left[ \left( a_{i,j}\right) _{\substack{ %
i\neq n-1  \\ j\neq n}}\right] =\underset{1\leq i,j\leq n-1}{\det }\left[
\det \left[ 
\begin{array}{cc}
a_{1,1} & a_{1,j+1} \\ 
a_{i+1,1} & a_{i+1,j+1}%
\end{array}%
\right] \right] _{\substack{ i\neq n-2  \\ j\neq n-1}}  \label{(2.7)}
\end{equation}

By using the formulas (\textrm{\ref{(2.4)})}-(\textrm{\ref{(2.7)})}, the
formula (\textrm{\ref{(2.3)})} will be as follows :

$\left\{ \left( a_{11}\right) ^{n-3}\right\} ^{2}\det \left[ \left(
a_{i,j}\right) _{1\leq i,j\leq n}\right] \det \left[ \left( a_{i,j}\right) 
_{\substack{ i\neq n-1,n  \\ j\neq n-1,n}}\right] =$

$\det \left[ 
\begin{array}{ll}
\underset{1\leq i,j\leq n-2}{\det }\left[ \left\vert 
\begin{array}{cc}
a_{1,1} & a_{1,j+1} \\ 
a_{i+1,1} & a_{i+1,j+1}%
\end{array}%
\right\vert \right] & \underset{1\leq i,j\leq n-1}{\det }\left[ \left\vert 
\begin{array}{cc}
a_{1,1} & a_{1,j+1} \\ 
a_{i+1,1} & a_{i+1,j+1}%
\end{array}%
\right\vert \right] _{\substack{ i\neq n-1  \\ j\neq n-2}} \\ 
\underset{1\leq i,j\leq n-1}{\det }\left[ \left\vert 
\begin{array}{cc}
a_{1,1} & a_{1,j+1} \\ 
a_{i+1,1} & a_{i+1,j+1}%
\end{array}%
\right\vert \right] _{\substack{ i\neq n-2  \\ j\neq n-1}} & \underset{1\leq
i,j\leq n-1}{\det }\left[ \left\vert 
\begin{array}{cc}
a_{1,1} & a_{1,j+1} \\ 
a_{i+1,1} & a_{i+1,j+1}%
\end{array}%
\right\vert \right] _{\substack{ i\neq n-2  \\ j\neq n-2}}%
\end{array}%
\right] .$

To simplify the above formula, we put this notation $d_{i,j}=\left\vert 
\begin{array}{cc}
a_{1,1} & a_{1,j+1} \\ 
a_{i+1,1} & a_{i+1,j+1}%
\end{array}%
\right\vert $, thus it will be as follows:%
\begin{equation*}
\left. 
\begin{array}{c}
\left\{ \left( a_{11}\right) ^{n-3}\right\} ^{2}\det \left[ \left(
a_{i,j}\right) _{1\leq i,j\leq n}\right] \det \left[ \left( a_{i,j}\right) 
_{\substack{ i\neq n-1,n  \\ j\neq n-1,n}}\right] = \\ 
\det \left[ 
\begin{array}{cc}
\underset{1\leq i,j\leq n-2}{\det }\left[ \left( d_{i,j}\right) \right] & 
\underset{1\leq i,j\leq n-1}{\det }\left[ \left( d_{i,j}\right) _{\substack{ %
i\neq n-1  \\ j\neq n-2}}\right] \smallskip \\ 
\underset{1\leq i,j\leq n-1}{\det }\left[ \left( d_{i,j}\right) _{\substack{ %
i\neq n-2  \\ j\neq n-1}}\right] & \underset{1\leq i,j\leq n-1}{\det }\left[
\left( d_{i,j}\right) _{\substack{ i\neq n-2  \\ j\neq n-2}}\right]%
\end{array}%
\right] \text{.}%
\end{array}%
\right.
\end{equation*}

We can write it as follows:%
\begin{eqnarray}
&&\left\{ \left( a_{11}\right) ^{n-3}\right\} ^{2}\det \left[ \left(
a_{i,j}\right) _{1\leq i,j\leq n}\right] \det \left[ \left( a_{i,j}\right) 
_{\substack{ i\neq n-1,n  \\ j\neq n-1,n}}\right]  \notag \\
&=&\det \left[ 
\begin{array}{cc}
\underset{1\leq i,j\leq n-1}{\det }\left[ \left( d_{i,j}\right) _{\substack{ %
i\neq n-1  \\ j\neq n-1}}\right] & \underset{1\leq i,j\leq n-1}{\det }\left[
\left( d_{i,j}\right) _{\substack{ i\neq n-1  \\ j\neq n-2}}\right] \\ 
\underset{1\leq i,j\leq n-1}{\det }\left[ \left( d_{i,j}\right) _{\substack{ %
i\neq n-2  \\ j\neq n-1}}\right] & \underset{1\leq i,j\leq n-1}{\det }\left[
\left( d_{i,j}\right) _{\substack{ i\neq n-2  \\ j\neq n-2}}\right]%
\end{array}%
\right] .  \label{(2.8)}
\end{eqnarray}

But by applying formula (\textrm{\ref{(2.3)})} for $\left( n-1\right) $ on $%
\left( d_{i,j}\right) _{1\leq i,j\leq n-1}$we get :%
\begin{equation*}
\left. 
\begin{array}{c}
\det \left[ \left( d_{i,j}\right) _{1\leq i,j\leq n-3}\right] \det \left[
\left( d_{i,j}\right) _{1\leq i,j\leq n-1}\right] = \\ 
\det \left[ 
\begin{array}{cc}
\underset{1\leq i,j\leq n-1}{\det }\left[ \left( d_{i,j}\right) _{\substack{ %
i\neq n-1  \\ j\neq n-1}}\right] & \underset{1\leq i,j\leq n-1}{\det }\left[
\left( d_{i,j}\right) _{\substack{ i\neq n-1  \\ j\neq n-2}}\right]
\smallskip \\ 
\underset{1\leq i,j\leq n-1}{\det }\left[ \left( d_{i,j}\right) _{\substack{ %
i\neq n-2  \\ j\neq n-1}}\right] & \underset{1\leq i,j\leq n-1}{\det }\left[
\left( d_{i,j}\right) _{\substack{ i\neq n-2  \\ j\neq n-2}}\right]%
\end{array}%
\right] .%
\end{array}%
\right.
\end{equation*}

Thus formula (\textrm{\ref{(2.8)})} will be as follaws:%
\begin{equation}
\left. 
\begin{array}{l}
\left\{ \left( a_{11}\right) ^{n-3}\right\} ^{2}\det \left[ \left(
a_{i,j}\right) _{1\leq i,j\leq n}\right] \det \left[ \left( a_{i,j}\right) 
_{\substack{ i\neq n-1,n  \\ j\neq n-1,n}}\right] = \\ 
\det \left[ \left( d_{i,j}\right) _{1\leq i,j\leq n-3}\right] \det \left[
\left( d_{i,j}\right) _{1\leq i,j\leq n-1}\right] .%
\end{array}%
\right.  \label{(2.9)}
\end{equation}

But%
\begin{equation*}
\det \left[ \left( a_{i,j}\right) _{\substack{ i\neq n-1,n  \\ j\neq n-1,n}}%
\right] =\det \left[ \left( a_{i,j}\right) _{1\leq i,j\leq n-2}\right] .
\end{equation*}

By applying formula (\textrm{\ref{(2.2)})} for $\left( n-2\right) $ on $%
\left( a_{i,j}\right) _{1\leq i,j\leq n-2}$, we get:%
\begin{equation*}
\left( a_{11}\right) ^{n-4}\underset{1\leq i,j\leq n}{\det }\left[ \left(
a_{i,j}\right) _{1\leq i,j\leq n-2}\right] =\underset{1\leq i,j\leq n-3}{%
\det }\left[ \det \left[ 
\begin{array}{cc}
a_{1,1} & a_{1,j+1} \\ 
a_{i+1,1} & a_{i+1,j+1}%
\end{array}%
\right] \right]
\end{equation*}%
the same as:%
\begin{equation}
\left( a_{11}\right) ^{n-4}\underset{1\leq i,j\leq n}{\det }\left[ \left(
a_{i,j}\right) _{1\leq i,j\leq n-2}\right] =\det \left[ \left(
d_{i,j}\right) _{1\leq i,j\leq n-3}\right] .  \label{(2.10)}
\end{equation}

Finally, formula (\textrm{\ref{(2.9)})} will be as follows:%
\begin{equation*}
\left\{ \left( a_{11}\right) ^{n-3}\right\} ^{2}\det \left[ \left(
a_{i,j}\right) _{1\leq i,j\leq n}\right] \det \left[ \left( a_{i,j}\right)
_{1\leq i,j\leq n-2}\right] =\det \left[ \left( d_{i,j}\right) _{1\leq
i,j\leq n-3}\right] \det \left[ \left( d_{i,j}\right) _{1\leq i,j\leq n-1}%
\right] .
\end{equation*}

And by using formula (\textrm{\ref{(2.10)})}, it will be as follows:%
\begin{equation*}
\left. 
\begin{array}{c}
\left\{ \left( a_{11}\right) ^{n-3}\right\} ^{2}\det \left[ \left(
a_{i,j}\right) _{1\leq i,j\leq n}\right] \det \left[ \left( a_{i,j}\right)
_{1\leq i,j\leq n-2}\right] = \\ 
\left( a_{11}\right) ^{n-4}\underset{1\leq i,j\leq n}{\det }\left[ \left(
a_{i,j}\right) _{1\leq i,j\leq n-2}\right] \det \left[ \left( d_{i,j}\right)
_{1\leq i,j\leq n-1}\right] .%
\end{array}%
\right.
\end{equation*}

At the end, we get:%
\begin{equation*}
\left( a_{11}\right) ^{n-2}\det \left[ \left( a_{i,j}\right) _{1\leq i,j\leq
n}\right] =\det \left[ \left( d_{i,j}\right) _{1\leq i,j\leq n-1}\right] .
\end{equation*}

And like this we have proved formula (\textrm{\ref{(2.2)})} for $n$.
\end{proof}

We can generalize theorem (\textrm{\ref{thor1})} by the following theorem:

\begin{theorem}
Let the $n$ square matrix $A=\left( a_{\left( i,j\right) }\right) _{1\leq
i,j\leq n}.$\newline
For all $n>2$ , we can generalize the formula (\textrm{\ref{(2.2)}) }as
follows:%
\begin{equation}
\left( a_{k,l}\right) ^{n-2}\det \left[ \left( a_{i,j}\right) _{1\leq
i,j\leq n}\right] =\underset{1\leq i,j\leq n-1}{\det }\left[ \det \left(
A_{i,j}\right) \right] \text{, \ \ \ \ }1\leq k,\text{ }l\leq n
\label{(2.11)}
\end{equation}%
when%
\begin{equation}
A_{\left( i,j\right) }=\left\{ 
\begin{array}{l}
\smallskip \left( 
\begin{array}{cc}
a_{i,j} & a_{i,l} \\ 
a_{k,j} & a_{k,l}%
\end{array}%
\right) \text{ \ \ \ \ \ \ \ \ if\ \ }j<l,i<k \\ 
\smallskip \left( 
\begin{array}{cc}
a_{i,l} & a_{i,j+1} \\ 
a_{k,l} & a_{k,j+1}%
\end{array}%
\right) \text{ \ \ \ \ \ if \ }j\geq l,i<k \\ 
\smallskip \left( 
\begin{array}{cc}
a_{k,j} & a_{k,l} \\ 
a_{i+1,j} & a_{i+1,l}%
\end{array}%
\right) \text{ \ \ \ \ if \ }j<l,i\geq k \\ 
\left( 
\begin{array}{cc}
a_{k,l} & a_{k,j+1} \\ 
a_{i+1,l} & a_{i+1,j+1}%
\end{array}%
\right) \text{ \ \ if \ }j\geq l,i\geq k%
\end{array}%
\right. \text{.}  \label{(2.12)}
\end{equation}
\end{theorem}

\begin{proof}
To prove formula (\textrm{\ref{(2.11)})}, we move the element $a_{k,l}$ from
its position to the position of the element $a_{1,1}$ in matrix $A$ by using
determinants properties in order to apply formula (\textrm{\ref{(2.2)}).}

So, we replace row $k$ and row $\left( k-1\right) $ by each other.Then, the
new row $\left( k-1\right) $ and row $\left( k-2\right) $ by each other, and
so on till row $k$ in matrix $A$ will be the first row. On the other side,
we replace column $l$ and column $\left( l-1\right) $ by each other.Then,
the new column $\left( l-1\right) $ and column $\left( l-2\right) $ by each
other, and so on till column $l$ in matrix $A$ will be the first column.We
get a new matrix $B$ that realises::%
\begin{equation}
\det A=\left( -1\right) ^{\left( k-1\right) +\left( l-1\right) }\det B.
\label{(2.13)}
\end{equation}

We apply formula (\textrm{\ref{(2.2)})} on matrix $B$, we get:%
\begin{equation}
\left( a_{k,l}\right) ^{n-2}\det B=\underset{1\leq i,j\leq n-1}{\det }\left[
\det \left( B_{i,j}\right) \right]  \label{(2.14)}
\end{equation}%
when:%
\begin{equation*}
B_{i,j}=\left\{ 
\begin{array}{l}
\smallskip \left( 
\begin{array}{cc}
a_{k,l} & a_{k,j} \\ 
a_{i,l} & a_{i,j}%
\end{array}%
\right) \text{ \ \ \ \ \ \ \ \ if\ \ }j<l,i<k \\ 
\smallskip \left( 
\begin{array}{cc}
a_{k,l} & a_{k,j+1} \\ 
a_{i,l} & a_{i,j+1}%
\end{array}%
\right) \text{ \ \ \ \ \ if \ }j\geq l,i<k \\ 
\smallskip \left( 
\begin{array}{cc}
a_{k,l} & a_{k,j} \\ 
a_{i+1,l} & a_{i+1,j}%
\end{array}%
\right) \text{ \ \ \ \ if \ }j<l,i\geq k \\ 
\left( 
\begin{array}{cc}
a_{k,l} & a_{k,j+1} \\ 
a_{i+1,l} & a_{i+1,j+1}%
\end{array}%
\right) \text{ \ \ if \ }j\geq l,i\geq k%
\end{array}%
\right. .
\end{equation*}%
by using determinant properties, we get:

\begin{equation*}
\det \left( B_{i,j}\right) =\left\{ 
\begin{array}{l}
\smallskip \left( -1\right) ^{2}\left\vert 
\begin{array}{cc}
a_{i,j} & a_{i,l} \\ 
a_{k,j} & a_{k,l}%
\end{array}%
\right\vert \text{ \ \ \ \ \ \ \ \ if\ \ }j<l,i<k \\ 
\smallskip \left( -1\right) \left\vert 
\begin{array}{cc}
a_{i,l} & a_{i,j+1} \\ 
a_{k,l} & a_{k,j+1}%
\end{array}%
\right\vert \text{ \ \ \ \ \ if \ }j\geq l,i<k \\ 
\smallskip \left( -1\right) \left\vert 
\begin{array}{cc}
a_{k,j} & a_{k,l} \\ 
a_{i+1,j} & a_{i+1,l}%
\end{array}%
\right\vert \text{ \ \ \ \ if \ }j<l,i\geq k \\ 
\left\vert 
\begin{array}{cc}
a_{k,l} & a_{k,j+1} \\ 
a_{i+1,l} & a_{i+1,j+1}%
\end{array}%
\right\vert \text{ \ \ if \ }j\geq l,i\geq k%
\end{array}%
\right. \text{.}
\end{equation*}

Thus, we write $\left( \det \left( B_{i,j}\right) \right) _{1\leq i,j\leq
n-1}$ as abock matrix%
\begin{equation*}
\left. 
\begin{array}{l}
\left( \det \left( B_{i,j}\right) \right) _{1\leq i,j\leq n-1}= \\ 
\left( 
\begin{array}{ll}
\left( k-1\right) \overset{\left( l-1\right) }{\overbrace{\left\{ \left(
\left( -1\right) ^{2}\left\vert 
\begin{array}{cc}
a_{i,j} & a_{i,l} \\ 
a_{k,j} & a_{k,l}%
\end{array}%
\right\vert \right) _{j<l,i<k}\right. }} & \overset{\left( n-l\right) }{%
\overbrace{\left( \left( -1\right) \left\vert 
\begin{array}{cc}
a_{i,l} & a_{i,j+1} \\ 
a_{k,l} & a_{k,j+1}%
\end{array}%
\right\vert \right) _{j>l,i<k}}} \\ 
\left( n-k\right) \left\{ \left( \left( -1\right) \left\vert 
\begin{array}{cc}
a_{k,j} & a_{k,l} \\ 
a_{i+1,j} & a_{i+1,l}%
\end{array}%
\right\vert \right) _{j<l,i>k}\right. & \left( \left\vert 
\begin{array}{cc}
a_{k,l} & a_{k,j+1} \\ 
a_{i+1,l} & a_{i+1,j+1}%
\end{array}%
\right\vert \right) _{j>l,i>k}%
\end{array}%
\right)%
\end{array}%
\right.
\end{equation*}%
and it can be written:%
\begin{equation*}
\left. 
\begin{array}{l}
\left( \det \left( B_{i,j}\right) \right) _{1\leq i,j\leq n-1}= \\ 
\left( -1\right) ^{\left( k-1\right) +\left( l-1\right) }\left( 
\begin{array}{ll}
\left( \left\vert 
\begin{array}{cc}
a_{i,j} & a_{i,l} \\ 
a_{k,j} & a_{k,l}%
\end{array}%
\right\vert \right) _{j<l,i<k} & \left( \left\vert 
\begin{array}{cc}
a_{i,l} & a_{i,j+1} \\ 
a_{k,l} & a_{k,j+1}%
\end{array}%
\right\vert \right) _{j>l,i<k} \\ 
\left( \left\vert 
\begin{array}{cc}
a_{k,j} & a_{k,l} \\ 
a_{i+1,j} & a_{i+1,l}%
\end{array}%
\right\vert \right) _{j<l,i>k} & \left( \left\vert 
\begin{array}{cc}
a_{k,l} & a_{k,j+1} \\ 
a_{i+1,l} & a_{i+1,j+1}%
\end{array}%
\right\vert \right) _{j>l,i>k}%
\end{array}%
\right)%
\end{array}%
\right.
\end{equation*}

According to formul (\textrm{\ref{(2.12)})}, we can write the last formula
as follows:%
\begin{equation*}
\left( \det \left( B_{i,j}\right) \right) _{1\leq i,j\leq n-1}=\left(
-1\right) ^{\left( k-1\right) +\left( l-1\right) }\underset{1\leq i,j\leq n-1%
}{\det }\left[ \det \left( A_{i,j}\right) \right] \text{.}
\end{equation*}

Thus, formula (\textrm{\ref{(2.14)})} will be:%
\begin{equation*}
\left( a_{k,l}\right) ^{n-2}\det B=\left( -1\right) ^{\left( k-1\right)
+\left( l-1\right) }\underset{1\leq i,j\leq n-1}{\det }\left[ \det \left(
A_{i,j}\right) \right] \text{.}
\end{equation*}

So%
\begin{equation*}
\left( a_{k,l}\right) ^{n-2}\left[ \left( -1\right) ^{\left( k-1\right)
+\left( l-1\right) }\det B\right] =\underset{1\leq i,j\leq n-1}{\det }\left[
\det \left( A_{i,j}\right) \right] \text{.}
\end{equation*}

By using formula (\textrm{\ref{(2.13)})}, we get formula (\textrm{\ref%
{(2.11)})}.
\end{proof}

\section{Application}

In this section, we'll show the main results we have found. These results
are the construction of an easy and simplified algorithm which compute the
determinant of any matrix (see \cite{Dodgson}). For this, we give the
following proposition:

\begin{proposition}
For $l=1,...,n:$%
\begin{equation}
\left( a_{\left( 1,l\right) }\right) ^{n-2}\det \left[ \left( a_{\left(
i,j\right) }\right) _{1\leq i,j\leq n}\right] =\underset{1\leq i,j\leq n-1}{%
\det }\left( \left\vert A_{i,j}\right\vert \right)  \label{(3.1)}
\end{equation}%
when%
\begin{equation*}
A_{i,j}=\left\{ 
\begin{array}{l}
\medskip \left( 
\begin{array}{cc}
a_{1,l} & a_{1,j+1} \\ 
a_{i+1,l} & a_{i+1,j+1}%
\end{array}%
\right) \text{ \ \ if \ }l\leq j \\ 
\left( 
\begin{array}{cc}
a_{1,j} & a_{1,l} \\ 
a_{i+1,j} & a_{i+1,l}%
\end{array}%
\right) \text{ \ \ \ \ if \ }j<l%
\end{array}%
\right. \text{.}
\end{equation*}
\end{proposition}

\begin{proof}
By putting $k=1$ in formula (\textrm{\ref{(2.11)})} we get formula (\textrm{%
\ref{(3.1)})}.
\end{proof}

\begin{algorithm}
This algorithm can be described in the following steps:\newline
$(01)$ $\ \ \ \ \ \ \ \ \ $Let $n\times n$ square matrix $A$ (we wish to
compute its determinant). \newline
$(02)$ \ \ \ \ \ \ \ If $n=2$, we compute $\left\vert A\right\vert $ by the
known method, else, \newline
$(03)$ \ \ \ \ \ \ \ if all the elements of the first row of matrix $A$ are
nil, then $\left\vert A\right\vert =0$, else, \newline
$(04)$ $\ \ \ \ \ \ \ $the first non nil elements in the first row is in $%
l^{th}$column. we form $\left( n-1\right) \times \left( n-1\right) $ square
matrix $B=\left( b_{i,j}\right) _{_{i,j=1}}^{n-1}$. Its elements are the
determinants of $2\times 2$ square matrix 
\begin{equation*}
b_{i,j}=\left\{ 
\begin{array}{l}
\medskip \left\vert 
\begin{array}{cc}
a_{1,l} & a_{1,j+1} \\ 
a_{i+1,l} & a_{i+1,j+1}%
\end{array}%
\right\vert \text{ \ \ if \ }l\leq j \\ 
\left\vert 
\begin{array}{cc}
a_{1,j} & a_{1,l} \\ 
a_{i+1,j} & a_{i+1,l}%
\end{array}%
\right\vert \text{ \ \ \ \ if \ }j<l%
\end{array}%
\right. \text{.}
\end{equation*}%
\newline
$(05)$ \ \ \ \ \ \ \ So,%
\begin{equation*}
\det A=\frac{\det B}{\left( a_{1,l}\right) ^{n-2}}
\end{equation*}%
\newline
$(06)$ \ \ \ \ \ \ \ Let $A=B$. We repeat the previous steps until we find
the determinant.\newline
\end{algorithm}

\end{document}